
\baselineskip=14pt
\parskip=10pt
\def\Tilde{\char126\relax}
\def\halmos{\hbox{\vrule height0.15cm width0.01cm\vbox{\hrule height
 0.01cm width0.2cm \vskip0.15cm \hrule height 0.01cm width0.2cm}\vrule
 height0.15cm width 0.01cm}}
\font\eightrm=cmr8  
\font\eighttt=cmtt8
\magnification=\magstephalf

\parindent=0pt
\overfullrule=0in
\bf
\centerline
{Dodgson's Determinant-Evaluation Rule Proved by }
\centerline
{TWO-TIMING MEN and WOMEN}
\bigskip
\centerline{ {\it Doron ZEILBERGER}\footnote{$^1$}
{\eightrm  \raggedright
Department of Mathematics, Temple University,
Philadelphia, PA 19122, USA. 
{\eighttt zeilberg@math.temple.edu \break
http://www.math.temple.edu/\Tilde zeilberg \quad
 ftp://ftp.math.temple.edu/pub/zeilberg} \quad .
\break
Supported in part by the NSF. Version of May 8, 1996.
First Version: April 16, 1996.
Submitted to the Wilf (Electronic) Festschrifft.
\break
Thanks are due to Bill Gosper for several corrections.
} 
}
\bigskip
\quad\quad\quad\quad\quad\quad\quad\quad\quad\quad
\quad\quad\quad\quad\quad\quad\quad\quad\quad\quad
{\it Bijections are where it's at  ---Herb Wilf}
\medskip
\quad\quad\quad\quad\quad{\it Dedicated to Master Bijectionist
Herb Wilf, on finishing 13/24 of his life}
\bigskip
\rm
I will give a bijective proof of the Reverend Charles Lutwidge
{\bf Dodgson's Rule}([D]):
$$
\det \left [ (a_{i,j})_{ {1 \leq i \leq n } \atop {1 \leq j \leq n }}
\right ]
\cdot
\det \left [
(a_{i,j})_{ {2 \leq i \leq n-1 } \atop {2 \leq j \leq n-1 }}
\right ] \, =
$$
$$
\det \left [
 (a_{i,j})_{ {1 \leq i \leq n-1 } \atop {1 \leq j \leq n-1 }}
\right ]
\cdot
\det \left [
 (a_{i,j})_{ {2 \leq i \leq n } \atop {2 \leq j \leq n }}
\right ]
\,-\,
\det \left [ (a_{i,j})_{ {1 \leq i \leq n-1 } \atop {2 \leq j \leq n }}
\right ]
\cdot
\det \left [
 (a_{i,j})_{ {2 \leq i \leq n } \atop {1 \leq j \leq n-1 }} 
\right ]\quad .
\eqno(Alice)
$$

Consider $n$ men, $1,2, \dots , n$,  and
$n$ women $1',2' \dots , n'$,
each of whom is married to exactly one member of the 
opposite sex. 
For each of the $n!$ possible (perfect) matchings
$\pi$, let
$$
weight(\pi):=sign(\pi) \prod_{i=1}^{n} a_{i,\pi(i)} \quad ,
$$
where $sign(\pi)$ is the sign of the corresponding permutation,
and for $i=1, \dots ,n$, Mr. $i$ is married to Ms. $\pi(i)'$.
 
Except for Mr. $1$, Mr. $n$, 
Ms. $1'$ and Ms. $n'$ all the persons have affairs.
Assume that each of the men in $\{2, \dots , n-1\}$ has exactly one
mistress amongst $\{2', \dots , (n-1)' \}$
and 
each of the women in $\{2', \dots , (n-1)'\}$ has exactly one
lover amongst $\{2, \dots , n-1 \}$\footnote{$^2$}
{\eightrm  \raggedright
Somewhat unrealistically,
a man's wife may also be his mistress, and equivalently,
a woman's husband may also be her lover.
}. 
For each of the $(n-2)!$ possible (perfect) matchings
$\sigma$, let
$$
weight(\sigma):=sign(\sigma) \prod_{i=2}^{n-1} a_{i,\sigma(i)} \quad ,
$$
where $sign(\sigma)$ is the sign of the corresponding permutation,
and for $i=2, \dots ,n-1$, Mr. $i$ is the lover of Ms. $\sigma(i)'$.
 
Let $A(n)$ be the set of all pairs $[ \pi , \sigma]$ as above, and let
$weight([\pi,\sigma]):=weight(\pi)weight(\sigma)$. 
The left side of $(Alice)$ is the sum of all the weights of the
elements of $A(n)$.
 
Let $B(n)$ be the set of pairs $[\pi,\sigma]$, where now
$n$ and $n'$ are unmarried but have affairs, i.e. $\pi$ is
a matching of $\{ 1, \dots , n-1\}$ to 
$\{ 1', \dots , (n-1)'\}$, and $\sigma$ is a matching of
$\{ 2, \dots , n\}$ to $\{ 2', \dots , n'\}$, and define the
weight similarly.
 
Let $C(n)$ be the set of pairs $[\pi,\sigma]$, where now
$n$ and $1'$ are unmarried and $1$ and $n'$ don't have affairs. 
i.e. $\pi$ is
a matching of $\{ 1, \dots , n-1\}$ to 
$\{ 2', \dots , n'\}$, and $\sigma$ is a matching of
$\{ 2, \dots , n\}$ to $\{ 1', \dots , (n-1)'\}$, and {\it now} define 
$weight([\pi,\sigma]):= -weight(\pi) weight(\sigma)$.
 
The right side of $(Alice)$ is the sum of all the weights of the
elements of $B(n) \cup C(n)$.
 
Define a mapping
$$
T: A(n) \rightarrow B(n) \cup C(n) \quad ,
$$
as follows. Given $[\pi,\sigma] \in A(n)$, define
an alternating sequence of men and women: 
$m_1:=n,w_1, m_2,w_2, \dots, m_r, w_r=1'$ or $n'$, such that
$w_i:=$wife of$( m_i)$, and $m_{i+1}:=$lover of$(w_i)$. This
sequence terminates, for some $r$, at either $w_r=1'$, or
$w_r=n'$, since then $m_{r+1}$ is undefined, as $1'$ and
$n'$ are lovers-less women. To perform $T$,
change the relationships 
$(m_1,w_1), (m_2,w_2), \dots , (m_r,w_r)$ from marriages to affairs
(i.e. Mr. $m_i$ and Ms. $w_i$ get divorced and become lovers, $i=1, \dots ,r$),
and change the relationships $(m_2,w_1), (m_3,w_2), \dots , (m_r,w_{r-1})$
from affairs to marriages. If $w_r=1'$ then 
$T([\pi,\sigma]) \in C(n)$, while if
$w_r=n'$ then $T([\pi,\sigma]) \in B(n)$.
 
The mapping $T$ is weight-preserving. Except for the sign, this
is obvious, since all the relationships have been preserved, only the
nature of some of them changed. I leave it as a pleasant exercise
to verify that also the sign is preserved.
 
It is obvious that $T: A(n) \rightarrow B(n) \cup C(n) \,\,$ 
is one-to-one. If it were onto, we would be done.
Since it is not, we need one more paragraph.
 
Call a member of $B(n) \cup C(n)$ {\it bad}
if it is not in $T(A(n))$. I claim that the sum of all the
weights of the bad members of $B(n) \cup C(n)$ is zero.
This follows from the fact that there
is a natural bijection $S$, easily constructed by the readers,
between the bad members of $C(n)$ and those of $B(n)$, such that
$weight(S([\pi, \sigma]))= - weight([\pi, \sigma])$. Hence 
the weights of the bad members of $B(n)$ and $C(n)$ cancel
each other in pairs, contributing a total of zero to the
right side of $(Alice)$. \halmos
 
A small Maple package, {\it alice}, containing programs
implementing the mapping $T$, its inverse, and the mapping
$S$ from the bad members of $C(n)$ to those of $B(n)$, is
available from my Home Page
{\tt http://www.math.temple.edu/\Tilde zeilberg}.
 
{\bf Reference}
 
[D] C.L. Dodgson, {\it Condensation of Determinants}, Proceedings
of the Royal Society of London {\bf 15}(1866), 150-155.
 
\bye